\documentclass{amsart}
\usepackage{amsmath,amssymb,verbatim}

\theoremstyle{definition}

\theoremstyle{remark}

\numberwithin{equation}{section}

\newcommand{\disk}{\mathbb{D}}

\newcommand{\up}{\upsilon}
\newcommand{\no}{\noindent}
\newcommand{\nah}{\nonumber}
\newcommand{\nl}{\newline}
\newcommand{\ka}{\kappa}
\newcommand{\la}{\lambda}

\newcommand{\ds}{\displaystyle}
\newcommand{\oln}[1]{\overline{#1}}
\newcommand{\pdisk}{\disk\times\disk}
\newcommand{\dee}{\partial}
\newcommand{\dbar}{\oln\partial}

\begin{document}

\title[Range of Berezin Transform]{Range of Berezin Transform}


\author{ N. V. Rao}
\address{The University Of Toledo, College of Arts and Sciences,\nl
 Department of Mathematics, Mail Stop 942. Toledo, Ohio 43606-3390,
USA}
\email{rnagise@math.utoledo.edu}




\date{\today}


\begin{abstract}
Let $\ds dA=\frac{dxdy}\pi$ denote the normalized Lebesgue area measure on the unit disk $\disk$ and $u$, a summable function on $\disk$.
$$B(u)(z)=\int_\disk u(\zeta)\frac{(1-|z|^2)^2}{|1-\zeta\oln z|^4}dA(\zeta)$$ is called the
Berezin transform of $u$.
Ahern \cite{a} described all the possible triples $\{u,f,g\}$  for which
$$B(u)(z)=f(z)\oln g(z)$$ where both $f,g$ are holomorphic in $\disk$. This result was crucial in solving a version of the zero product problem for Toeplitz operators on the Bergman space.

The natural next question was to describe all functions in the range of Berezin Transform which are of the form $$\sum_{i=1}^Nf_i\oln g_i$$ where $f_i,g_i$ are all holomorphic in $\disk$. We shall give a complete description
of all such $u$ and the corresponding $f_i,g_i,1\leq i\leq N$.
Further we give very simple proof of the result of Ahern \cite{a} and the recent results of \v Cu\v ckovi\'c and Li \cite{bz} where
they tackle the special case when $N=2$ and $g_2=z^n$.

\end{abstract}

\maketitle

\section{Introduction}
Let $\ds dA=\frac{dxdy}\pi$ denote the normalized Lebesgue area measure on the unit disk $\disk$ in the
complex plane.
For any $u$, a summable function on $\disk$,
$$B(u)(z)=\int_\disk u(\zeta)\frac{(1-|z|^2)^2}{|1-\zeta\oln z|^4}dA(\zeta)\eqno(1)$$ is called the
Berezin transform of $u$.

\no{\bf Theorem A:} (Ahern \cite{a}) {\em If $u\in L^1(\disk)$ and
$$B(u)(z)=f(z)\oln g(z)$$
where both $f,g$ are holomorphic in $\disk$ and not constant, then there exists an automorphism
$\phi$ of $\disk$ and two polynomials $p,q$ each of degree at most 2 and the degree of $pq$ is at most 3
such that $f=p(\phi),g=q(\phi)$.}

By the standard trick of complexification  we obtain $B(u)(z,w)$ a holomorphic function of two complex variables defined in the bidisk $\pdisk$
as follows:
$$B(u)(z,w)= \int_\disk u(\zeta)\frac{(1-zw)^2}{(1-\oln\zeta z)^2(1-\zeta w)^2}dA(\zeta)$$ and notice that (1) can
be rewritten as
$$B(u)(z)=B(u)(z,\oln z).$$ Further we can write
$$B(u)(z,w)=\sum_{k=0}^\infty f_k(z)w^k$$ where $f_k$ is holomorphic in $\disk$ for every $k\geq 0$.
For any function $v(z)$, let
\begin{displaymath}
\begin{array}{ccl}
\dee v & = & \frac{\partial v}{\partial z}\\
\dbar v& =& \frac{\bar\partial v}{\partial\oln z}\\
\Delta&=&\dee\dbar,\mbox{ the Laplace Operator.}
\end{array}
\end{displaymath}
\vskip 10pt
\no {\em We say $B(u)(z)$ and also $B(u)(z,w)$ is of {\bf finite
rank} if the vector space generated by $f_k$ is finite dimensional.} The following remarks are
not difficult to check.
\vskip 10pt
\no{\bf Remark 1:} If $B(u)$ is of finite rank, then it can be written as
$$\sum_{i=1}^Nf_i\oln g_i,\mbox{ where }f_i,g_i\mbox{ are holomorphic in }\disk,$$
and conversely.
\vskip 10pt
\no{\bf Remark 2:} $$\sum_{i=1}^Nf_i\oln g_i$$ is of rank $N$ if and only if
the set of functions $\{f_1,f_2,\ldots f_N\}$ is linearly independent and so is
$\{g_1,g_2,\ldots g_N\}$.
\vskip 10pt
\no This terminology could be applied to any holomorphic function of two complex variables in the bidisk. For example the function $f(z)\oln g(z)$ is of
rank 1 because after complexification it will be $f(z)\oln g(w)=\sum_{k=0}^\infty c_kf(z)w^k$.
\vskip 10pt
\no {\em Ahern's theorem can be thought of as characterizing all functions of rank one in the range of Berezin transform on $\disk$.}
\vskip 10pt
\no{\em The theorem of \cite{bz}, including both cases, comes under the case of rank not exceeding 3}
because any harmonic function can be written as $f_1\oln g_1+f_2\oln g_2$ where $f_1,g_2$  are holomorphic and
$g_1=f_2=1$ and so case of (2) is covered under rank not exceeding 3 and case of (3) is covered under  rank not exceeding 2.
\section{Main Theorem}
\no{\bf  Theorem 1:} {\bf If $B(u)$ is of finite rank not exceeding $N$, then there exist finitely many points $a_i,1\leq i\leq N$ in $\disk$ such that $$u(\zeta)=h(\zeta)+\sum_{i=1}^ND_i\ln|\zeta-a_i|+\frac{E_i}{(\zeta-a_i)}+\frac{F_i}{\overline{(\zeta-a_i)}}$$
where  $D_i,E_i,F_i,1\leq i\leq N$ are constants, many and even all of them could vanish, and $h$ is a summable harmonic function.}

\no{\bf Proof.}
We follow Ahern by applying Laplacian to $B(u)$,
\begin{displaymath}
\begin{array}{rcl}
\Delta B(u)&=& \sum_{k=1}^\infty f_k'(z)k\oln z^{k-1}\\
 &=& \int_\disk u(\zeta)\Delta_z\frac{(1-|z|^2)^2}{|1-\zeta\oln z|^4}dA(\zeta)
\end{array}
\end{displaymath}  where
$\Delta_z$ denotes Laplacian with respect to $z$. As noted by Ahern\cite{a} there is a remarkable  symmetry for the
Berezin kernel
$$\Delta_z \frac{(1-|z|^2)^2}{|1-\zeta\oln z|^4}=\Delta_\zeta\frac{(1-|\zeta|^2)^2}{|1-\oln\zeta z|^4}$$
and so
\begin{eqnarray}
\sum_{k=1}^\infty f_k'(z)k\oln z^{k-1} &=& \int_\disk u(\zeta)\Delta_\zeta\frac{(1-|\zeta|^2)^2}{|1-\oln\zeta z|^4}dA(\zeta)\nah\\
&=&\int_\disk u(\zeta)\Delta_\zeta\frac{(1-|\zeta|^2)^2}{(1-\oln\zeta z)^2(1-\zeta \oln z)^2} dA(\zeta)\nah\\
&=&\int_\disk u(\zeta)\Delta_\zeta\frac{(1-|\zeta|^2)^2}{(1-\oln\zeta z)^2}
\sum_{k=0}^\infty(k+1)\zeta^k \oln z^k dA(\zeta)\nah\\
&=&\sum_{k=0}^\infty(k+1) \oln z^k\int_\disk u(\zeta)\Delta_\zeta\frac{(1-|\zeta|^2)^2}{(1-\oln\zeta z)^2}
\zeta^kdA(\zeta).\nah
\end{eqnarray} From this we get, after complexifying both sides
$$\sum_{k=1}^\infty f_k'(z)kw^{k-1}= \sum_{k=0}^\infty(k+1) w^k\int_\disk u(\zeta)\Delta_\zeta\frac{(1-|\zeta|^2)^2}{(1-\oln\zeta z)^2}
\zeta^kdA(\zeta)$$
and equating the coefficients of $w^k$,

$$f'_{k+1}(z)=\int_\disk u(\zeta)\Delta_\zeta\frac{(1-|\zeta|^2)^2}{(1-\oln\zeta z)^2}
\zeta^kdA(\zeta).\eqno(4)$$ Let us write the power series expansion for $f_{k+1}'$ as follows:
$$f'_{k+1}(z)=\sum_{l=0}^\infty a_{k,l}z^l.$$ From (4) we deduce
$$ f'_{k+1}(z)=\sum_{l=0}^\infty (l+1)z^l\int_\disk u(\zeta)\Delta_\zeta(1-|\zeta|^2)^2\oln\zeta^l
\zeta^kdA(\zeta)$$ and so
$$\frac{a_{k,l}}{l+1} = \int_\disk u(\zeta)\Delta_\zeta(1-|\zeta|^2)^2\oln\zeta^l
\zeta^kdA(\zeta){\rm\quad for \quad every\quad} k,l\geq 0.$$ We are given that the vector space generated by
$f_k$ is finite dimensional and so same goes for the vector space generated by $f_k'$ which means the
matrix $\{a_{k,l}\}$ is of finite rank and therefore the matrix $\left\{\ds\frac{a_{k,l}}{l+1}\right\}$ is also
of finite rank since the column space is unaltered. Now we see the distribution $F$ with support in $\oln\disk$ defined by
$$F(\psi)=\int_\disk u(\zeta)\Delta_\zeta(1-|\zeta|^2)^2\psi(\zeta)dA(\zeta)$$ for any $C^\infty$ function $\psi$
in the complex plane, satisfies the
Luecking \cite{l} condition: The matrix $\ds\left\{F(\zeta^l\oln\zeta^k)\right\}$ is of  rank,  $\leq N$.

By the theorem 3.1 of
Alexandrov and Rozenblum \cite{ar} we have that, such a distribution has only finite support, that is there exist at most $N$ points $a_i$ in  $\disk$ and linear differential operators $L_i$ such that
$$F=(1-|\zeta|^2)^2\Delta_\zeta u(\zeta)=\sum_{i=1}^NL_i(\delta(\zeta-a_i))$$ where $\delta$ is the Dirac delta
function.  This implies that the distribution $\Delta u\equiv 0$  except at $a_i$ and so $u$ is harmonic in $\disk$ with isolated singularities at $a_i$. Hence in a neighborhood of each $a_i$, it will have a Laurent series like $$c_0\ln|\zeta-a_i|+c_1/(\zeta-a_i)+c_2/\overline{(\zeta-a_i)}.$$ Higher powers will be absent since $u$ is summable. This proves the main theorem.

\section{Applications}
 The Berezin transforms of basic functions $\ln|\zeta|$, $\ds\frac1\zeta$, and
$\ds\frac1{\oln\zeta}$   are calculated in Ahern's paper\cite{a}. For example $B(\ln|\zeta|)=\ds\frac{z\oln z-1}2$ and
$\ds B(1/\zeta)=2\oln z-z\oln z^2$ and the rest can be computed by fractional linear transformations. For example it is an easy exercise to calculate Berezin transform of $\ln|\zeta-a|$ and $\ds\frac1{\zeta-a}$ because
 $$B(u(\phi_a))=B(u)(\phi_a)\eqno(5)$$ where
$$\phi_a(\zeta)=\ds\frac{\zeta-a}{1-\oln a\zeta}$$ and if $v$ is harmonic and summable,
$$B(v)=v.\eqno(6)$$
For example
$$B(\ln|\zeta-a|)=B(\ln|\phi_a(\zeta)|)+B(\ln|1-\oln a\zeta|)=\frac{\phi_a\oln\phi_a-1}2+
\ln|1-\oln a\zeta|,$$
$$2B(\ln|\zeta-a|-\ln|1-\oln a\zeta|)=\phi_a\oln\phi_a,\eqno(6a)$$ and
$$(1-a\oln a) B\left(\frac1{\zeta-a}\right)=B\left(\frac{1-\oln a\zeta}{\zeta -a}+\oln a\right)=B\left(\frac1{\phi_a}\right)+\oln a
=2\oln\phi_a-\phi_a\oln\phi_a^2+\oln a,$$
$$B\left(\oln a+2\oln\phi_a-\frac{1-a\oln a}{\zeta-a}\right)=\phi_a\oln\phi_a^2.\eqno(6b)$$
 So we conclude
 \vskip 10pt

 \no{\bf Theorem 2.} If $u$ is summable and $B(u)$ of finite rank $N$, then  $$B(u) = h+\sum_{i=1}^N\phi_{a_i}\oln\phi_{a_i}(D_i+E_i\phi_{a_i}+F_i\oln\phi_{a_i})\eqno(7)$$ where $h$ is a summable harmonic function and $a_i\in\disk$ and $D_i,E_i,F_i$ are constants.
 \vskip 10pt
 \no {\bf Corollary 1.} (\v Cu\v ckovi\'c and Li, \cite{bz}) If $u$ is summable and $B(u)$ is harmonic, then $u$ is harmonic.
\vskip 10pt

\no Proof. As noted this means $B(u)$ is of rank at most 2. So by  Theorem 2 there exist at most two points $a_1,a_2$ in $\disk$ such that
 $$B(u)=h+\sum_{i=1}^2\phi_{a_i}\oln\phi_{a_i}(A_i+B_i\phi_{a_i}+C_i\oln\phi_{a_i})$$ where $h$ is harmonic
 and summable. So $B(u)=v$ is summable and since it is given to be  harmonic, we have $B(v)=v$ and therefore $B(u-v)=0$. It is well known that $B$ is injective and hence $u=v$ and $u$ is harmonic. QED.
 \nl{\bf Comment.} Prof. Ahern communicated to me that this fact was noted by him long ago and his proof even simpler than what is presented here,  goes as follows: Main point is to prove that $B(u)$ is summable.
 Using Lemma 6.23 of [\cite{kz}, p 148] along with the notation there,  we have
 $$\tilde\Delta B(u)=8(B(u)-B_1(u)).\eqno(\rm afr)$$
 But harmonicity of $B(u)$ makes the LHS zero and so
 $B(u)=B_1(u)$. Since $B_1(u)$ is summable, $B(u)$ is summable. The equation (afr) was already noted in
 \cite{afr}.
 \vskip 10pt

\no Now we apply Theorem 2 to prove Theorem A:\vskip 10pt

\no{\bf Proof of Theorem A.}  So by the hypothesis of Theorem A, $B(u)$ has rank 1 and so from   Theorem 2 follows the existence of a point $a\in\disk$ such that $$f(\zeta)\oln g(\zeta)=h(\zeta)+(A\phi_a(\zeta)+B\phi_a^2(\zeta))\oln\phi_a(\zeta)+C\phi_a(\zeta)\oln\phi_a^2(\zeta)$$ where $h$ is a summable harmonic function
and $A,B,C$ are constants and\newline $\phi_a(\zeta)=\ds\frac{\zeta-a}{1-\oln a\zeta}$, an automorphism of $\disk$. By (5) and a change of variable $z=\phi_a(\zeta)$, we
get $$F(z)\oln G(z)=H(z)+(Az+Bz^2)\oln z+Cz\oln z^2\eqno(8)$$ where $F(z)=f(\zeta),G(z)=g(\zeta),H(z)=h(\zeta)$.
Since $H$ is harmonic in $\disk$, we can write it in a unique way as $$K(z)+\oln L(z)$$ where $K,L$ are holomorphic in $\disk$ and $L(0)=0$. Also let us write the Taylor series $$\oln G(z)=\sum_{k=0}^\infty g_k\oln z^k,\oln L(z)=\sum_{k=1}^\infty l_k\oln z^k.\eqno(9)$$ Now (8) can be written as
$$F(z)\oln G(z)\equiv\sum_{k=0}^\infty g_kF(z)\oln z^k\equiv K(z)+(l_1+Az+Bz^2)\oln z+(l_2+Cz)\oln z^2+\sum_{k=3}^\infty l_k\oln z^k.$$ Comparing the coefficient of $\oln z^k$ on both sides, we obtain
$$g_0F(z)=K(z),g_1F(z)=l_1+Az+Bz^2,g_2F(z)=l_2+Cz,g_kF(z)=l_k{\rm\quad for\quad}k\geq 3.$$
So if $g_k\neq 0$ for some $k>2$, $F$ and so $f$ will be constant.

\no Since $f$ is not constant,   we get
$g_k=0$ for all $k>2$ and that means $G$ is a polynomial of degree $\leq 2$. Now one of $g_1,g_2$ is
different from zero for otherwise $G$ would be constant.

\no Now if $g_2=0,g_1\neq 0$,
G is a polynomial of degree 1 and $ F(z)=(l_1+Az+Bz^2)/g_1$ is a polynomial of degree at most 2. On the other hand $g_2\neq 0$, $F(z)=(l_2+Cz)/g_2$ is a polynomial of degree 1 and $G$ is of degree 2.
This proves the theorem A by noting $$F(\phi_a(\zeta))=f(\zeta),\quad G(\phi_a(\zeta))=g(\zeta),$$
$F(z)=p(z),G(z)=q(z)$ the promised polynomials.
\vskip 10pt

 \no Let us end this section with a lemma that will be useful in the next section.\nl
 {\bf Lemma 1:} Let $\ds\phi(z)=\frac{z-a}{1-\oln az}, |a|<1$ and $A,B,C$ be constants, not all zero and
 $$\psi=(A\phi+B\phi^2)\oln\phi+C\phi\oln\phi^2.$$ There are two possibilities.
 \begin{itemize}
 \item[a)] $BC\neq 0$, and there exist functions $u_1,u_2$ both summable
 such that $$\psi=B(u_1)+B(u_2)\mbox{ and
} B(u_1)=f_1\oln g_1,B(u_2)=f_2\oln g_2$$ where $f_1,g_1,f_2,g_2$ are holomorphic everywhere except
 at $b=\ds\frac1{\oln a}$. Further $f_1,g_2$ have  a pole of order 2 and
 $f_2,g_1$ have a pole of order 1 at $b$.
 \item[b)] $BC=0$ and  there exists a summable function $u$ such that
 $$\psi=B(u)=f\oln g$$ where $f,g$ are holomorphic everywhere except at $b$ where each has a pole.

 \end{itemize}
\no{\bf Proof.} Let us assume $BC\neq 0$ which means $B\neq 0$ and $C\neq 0$. Now we set
$$f_1=A\phi+B\phi^2,g_1=\phi, f_2=C\phi, g_2=\phi^2.$$ Now from (6a) and (6b) follows
 the existence of  summable functions
$u_1$ and $u_2$ such that $$B(u_1)=f_1\oln g_1,B(u_2)=f_2\oln g_2.$$ Easy to check all of  a)
is true.

\no If $BC=0$, one of $B,C$ is zero. If $B=0,C\neq 0$, we set
$f=\phi, g=\oln A\phi+\oln C\phi^2$ and if $B\neq 0,C=0$, set
$f=A\phi+B\phi^2,g=\phi$. And if $B=C=0$, we set $f=A\phi, g=\phi$.
Now easy to check in each instance from (6a) and (6b) that there exists a smmable function
$u$ such that $B(u)=f\oln g$ and easy to check that all of b) is true.

\vskip 10pt

\section{The case of rank N, any positive integer.}
\no {\bf Theorem 3.} If $u\in L^1(\disk)$ and $B(u)$ is of rank $N$, then either $\ds u=\sum_{i=1}^Nu_i$ where each $u_i$ is summable and each of $B(u_i)$ is of rank one or $u_1+u_2$ is summable and harmonic and for $i>2$,
$u_i$ is summable and $B(u_i)$ is of rank one.
\nl [I think it is possible for a harmonic function to be summable but not its conjugate.]
\vskip 10pt
\no{\bf Proof.} Theorem 2 implies there exist at most $N$ distinct automorphisms
$$\phi_i=\ds\frac{z-a_i}{1-\oln a_iz}$$
of $\disk$ and a harmonic function $h$ and
constants $D_i,E_i,F_i, 1\leq i\leq$, allowing the possibility for a lot of these
constants to vanish, such that
$$B(u) = h + \sum_{i=1}^N(D_i\phi_i + E_i\phi_i^2)\oln\phi_i + F_i\phi_i\oln\phi_i^2 .\eqno(10)$$
Write $h$ uniquely as $K+\oln L$ where $K,L$ are holomorphic and $L(0)=0$.
\nl{\bf Fix an $i$.} From Lemma 1, there are three possibilities for the function
$$\psi_i=(D_i\phi_i + E_i\phi_i^2)\oln\phi_i + F_i\phi_i\oln\phi_i^2,$$
either it could be zero or there exist two summable functions $u,v$ such that
$$\psi_i=B(u)+B(v), B(u)=f\oln g, B(v)=p\oln q$$ such that $f,g,p,q$ are holomorphic everywhere except
at $\frac1{\oln a_i}=b_i$ where $f,q$ have a pole of order 2 and $g,p$ have a pole of order 1.
Or there exists just one summable function $u$ such that
$$\psi_i=B(u)=f\oln g$$ where $f,g$ are holomorphic everywhere except at $b_i$
where each of $f,g$ has a pole of order $\leq 2$.

\no Now from (10), we deduce the existence of finitely many summable functions $u_j,1\leq j\leq M$
such that $$B(u_j)=f_j\oln g_j$$ and
$$B(u)=h+\sum_{j=1}^MB(u_j)=h+\sum_{j=1}^Mf_j\oln g_j,\eqno(11)$$
{\bf where each of the functions $f_j,g_j$ is holomorphic in the entire plane except at a single point,where
each has a pole and
which varies with the index $j$ and if two pairs $(f_j,g_j)$ and $(f_k,g_k)$ have poles at
the same point, this can happen only for two pairs both coming from one $\psi_i$, the orders of
the pole differ from $f_j$ to $f_k$ and from $g_j$ to $g_k$.}
\vskip 10pt

 \no Summing up the above discussion,
we have $$B(u)=K+\oln L+\sum_{i=1}^MB(u_i),\eqno(12)$$
$$B(u_i)=f_i\oln g_i,\quad 1\leq i\leq M,$$ where each $f_i$ and $g_i$ is meromorphic in the entire
plane with a single pole, of course the pole may vary with $i$, and if two of the $f_i$ happen to have a pole at the same point, one of them will have pole of order 2 and the other of order 1. Same can be said about the $g_i$.These are the contributions from
the same $\psi_j$.\vskip 10pt

\no Let us write the Taylor series for $g_i,L$,
$$g_i(z)=\sum_{k=0}^\infty\oln\beta_{i,k}z^k,L(z)=\sum_{k=1}^\infty\oln l_kz^k.$$

\no Let $w_k$ denote the coefficient of $\oln z^k$ for the function $B(u)$. Then $N$ is the dimension of the vector space $W$ generated by $w_k$. We list
\begin{displaymath}
\begin{array}{rl}
w_0=&K+\sum_{i=1}^M\beta_{i,0}f_i\\
w_k=&l_k+\sum_{i=1}^M\beta_{i,k}f_i\mbox{\quad for }k\geq 1
\end{array}
\end{displaymath}using (10).

\no Evidently dimension of $W_1$
space generated by $w_k,k>0$ would be $\leq N$. We claim  $1,f_1,f_2,\ldots f_M$ are linearly independent for if not,
there would exist constants $\xi_0,\xi_1,\ldots\xi_M$ such that
$$\xi_0+\xi_1 f_1\ldots+\xi_m f_m\equiv 0.$$ As noted earlier each of $f_i$ is meromorphic in the plane with a single pole and the poles are either different or of different orders if they happen to be the same. Hence each $\xi_i=0$ for $i>0$. Hence $\xi_0=0$. \vskip 10pt

\no Therefore
$w_k$ is represented by the row $(l_k,\beta_{1,k},\beta_{2,k},\ldots\beta_{M,k})$. Let $W_1$ also denote the matrix whose rows are $w_k,k>0$. Then rank of $W_1$ is $\leq N$ and the matrix $C$ obtained by dropping the first column
of $W_1$ would have rank $\leq N$. But rank of $C$ is $M$. For otherwise there would exist constants
$\xi_1,\ldots\xi_M$, not all zero,such that
$$\sum_{i=1}^M\xi_i\beta_{i,k}=0,\forall k>0,$$
and so
$$\sum_{i=1}^M\xi_i(\oln g_i-\beta_{i,0})=0.$$
But again using the poles argument, we see that
all $\xi_i=0$. Hence the rank of $C$ is $M$ and $M\leq N$. Evidently
the rank of $B(u)=K+\oln L+\sum_{i=1}^MB(u_i)$ is less than or equal to
$1+1+M$. So $$N-2\leq M\leq N.$$ \vskip 10pt

\no If $M=N-2$, the theorem is proved because $B(u)=B(h)+\sum_{i=1}^MB(u_i)$ where $h$ is summable and harmonic.

\no Assume $M=N-1$. There are two possibilities depending on the rank of $W_1$. Its rank lies between
$M$ and $N$. If its rank is $M=N-1$, the first column will depend on the rest and so we see that
$\oln L$ is a linear combination of $\oln g_i-\oln g_i(0),1\leq i\leq M$. Therefore there exist constants
$\la_i$ such that
$$\oln L=\sum_{i=1}^M\la_i(\oln g_i-\oln g_i(0)),\eqno(13)$$ from which and (12), we obtain
$$B(u)=K+\sum_{i=1}^M\la_i(\oln g_i-\oln g_i(0))+\sum_{i=1}^Mf_i\oln g_i,$$
$$=K-\sum_{i=1}^M\la_i\oln g_i(0)+\sum_{i=1}^M(f_i+\la_i)(\oln g_i)$$
$$=K-\sum_{i=1}^M\la_i\oln g_i(0)+\sum_{i=1}^MB(u_i+\la_i\oln g_i).$$
Now by (13), $L$ is summable and since $K+\oln L$ is summable,\newline
 we have
$K-\sum_{i=1}^M\la_i\oln g_i(0)=K_1$ is summable, further since it is holomorhic,
$B(K_1)=K_1$. Therefore
$$B(u)=B(K_1)+\sum_{i=1}^MB(u_i+\la_i\oln g_i)$$ proving the Theorem in this case.

\no The remaining possibility is that the rank of $W_1$ is $N$. Then $w_0$ is a linear combination of
some of the rows of $W_1$, let us say
$$w_0=\sum_{i=1}^n\mu_iw_i.$$ Hence
$$w_0=K+\sum_{j=1}^M\beta_{j,0}f_j=\sum_{i=1}^n\mu_i\left(l_i+\sum_{j=1}^M\beta_{j,i}f_j\right),$$
$$K=\ka_0+\sum_{j=1}^M\ka_jf_j.\eqno(14)$$ Therefore from (12), we have
$$B(u)=\oln L+\ka_0+\sum_{j=1}^M\ka_jf_j+\sum_{j=1}^Mf_j\oln g_j$$
$$=\oln L+\ka_0+\sum_{j=1}^Mf_j(\oln g_j+\ka_j)$$
$$=\oln L+\ka_0+\sum_{j=1}^MB(u_j+\ka_j f_j).$$ Now by (14), $K$ is summable
and so $\oln L_1=\oln L+\ka_0$ is summable and further since it is anti-holomorphic,
$B(\oln L_1)=\oln L_1.$  Therefore
$$B(u)=B(\oln L_1)+\sum_{j=1}^MB(u_j+\ka_j f_j)$$ proving the Theorem in this case.
\vskip 10pt
\no Assume $M=N$. In this case first column of $W_1$ depends on the other columns which are columns of
$C$ and this means $\oln L$ is a linear combination of $\oln g_i-\oln g_i(0),1\leq i\leq M$.
Following the same argument as in (13), we obtain
$$B(u)=K-\sum_{i=1}^M\la_i\oln g_i(0)+\sum_{i=1}^M(f_i+\la_i)\oln g_i.$$ Now by introducing
new functions in place of the old $K,L,f_i,g_i$, we do not change $B(u)$ and nor the new
functions. They have the same properties as the old ones. So Let
$$\tilde K=K-\sum_{i=1}^M\la_i\oln g_i(0),\tilde L=0,\tilde f_i=f_i+\la_i,\tilde g_i=g_i.$$
Here $$\tilde u_i=u_i-\la_i\oln g_i$$
Now the coefficient $w_k$ of $\oln z^k$ of $B(u)$ would look like
\begin{displaymath}
\begin{array}{rl}
w_0=&\tilde K+\sum_{i=1}^M\beta_{i,0}\tilde f_i\\
w_k=&\sum_{i=1}^M\beta_{i,k}\tilde f_i\mbox{ for }k\geq 1.
\end{array}
\end{displaymath}
Since the rank of $B(u)$ is equal to the rank of the matrix $(\beta_{i,k})=C$, as argued in the
above paragraph $w_0$ is a linear combination of $w_k,0<k\leq n$ for some $n$. Therefore
we get
$$\tilde K+\sum_{i=1}^M\beta_{i,0}\tilde f_i=\sum_{k=1}^n\mu_i\sum_{i=1}^M\beta_{i,k}\tilde f_i$$
from which we get
$$\tilde K=\sum_{i=1}^M\ka_i\tilde f_i,$$
and $$ B(u)=\sum_{i=1}^M\ka_i\tilde f_i + \sum_{i=1}^M\tilde f_i\oln g_i$$
$$=\sum_{i=1}^M\tilde f_i(\oln g_i+\ka_i)$$
$$=\sum_{i=1}^MB(\tilde u_i+\ka_i\tilde f_i).$$
This proves the Theorem 3 completely.

\end{document}